\documentclass{article}
\newtheorem{conjecture}{Conjecture}

\newtheorem{lemma}{Lemma}
\newtheorem{algorithm}{Algorithm}

\usepackage{amssymb}
\usepackage{array}
\usepackage{float}
\usepackage{url}

\usepackage{array}
\newlength{\cellwid}
\makeatletter
\newenvironment{latinsq}[1][00]{%
 \def\centcol{\centering \let\\=\tabularnewline \CellStrut}
 \arraycolsep0pt
 \setbox\@tempboxa\hbox{#1}%
 \cellwid\ht\@tempboxa  \advance\cellwid\dp\strutbox
 \ifdim\cellwid<\wd\@tempboxa 
   \@tempdima .5\wd\@tempboxa \advance\@tempdima -.5\cellwid 
   \advance\cellwid \@tempdima \advance\@tempdima\dp\strutbox
 \else
   \@tempdima\dp\strutbox
 \fi
 \edef\CellStrut{\vrule
   width\z@ height\the\cellwid depth\the\@tempdima \relax}
 \advance\cellwid\@tempdima \advance\cellwid-2\arraycolsep
 \array{|>{\centcol}p{\cellwid}|*{20}{>{\centcol}p{\cellwid}|}}}%
 {\endarray}
\makeatother

\title{The size of the smallest uniquely completable set in order
8 Latin squares}
\author{Richard Bean}
\date{}

\begin{document}
\maketitle
\begin{center}
Institute for Studies in Theoretical Physics and Mathematics \\
Department of Mathematics \\
PO Box 19395-5746 \\
Tehran, I.R. Iran \\
rwb at ipm-dot-ir \\
\end{center}

\begin{abstract}
In 1990, Kolesova, Lam and Thiel determined the 283,657 main
classes of Latin squares of order 8.  Using techniques to
determine relevant Latin trades and integer programming, we
examine representatives of each of these main classes and
determine that none can contain a uniquely completable set of size
less than 16. In three of these main classes, the use of trades
which contain less than or equal to three rows, columns, or
elements does not suffice to determine this fact.  We closely
examine properties of representatives of these three main classes.
Writing the main result in Nelder's notation for critical sets, we
prove that scs(8)=16.
\end{abstract}

\section{Introduction}
A {\it Latin square} $L$ of order $n$ is an $n \times n$ array of
entries $\{ (i,j;k) \}$ such that each row and column of $L$
contains each of $n$ possible elements exactly once.  We will
refer to the Latin square based on the group table of
$\mathbb{Z}_n$ simply as $\mathbb{Z}_n$.  A {\it uniquely
completable (UC)} set $U$ is a subset of a Latin square $L$ such
that $L$ is the only superset of $C$ which is a Latin square. A
{\it critical set} $C$ of $L$ is a subset of $L$ such that $C$ is
uniquely completable and no subset of $C$ has this property.
Critical sets were introduced by Nelder \cite{nelder}.

A related concept is that of the {\it Latin trade}.  In this
paper, we will consider a Latin trade to be the set of entries in
which two Latin squares of the same order differ.  An {\it
intercalate} is a Latin trade of size 4.  The connection between
uniquely completable sets and Latin trades is well-known, for
example see \cite{aka}, \cite{def}.  It is expressed in the
following lemma.

\begin{lemma}
In any uniquely completable set $U$ for a Latin square $L$, each
trade in $L$ must intersect $U$ in at least one position.
\end{lemma}

scs($n$) is a function defined by Nelder as the smallest size of a
critical set in any $n \times n$ Latin square.  Nelder determined
that scs(3)=2, Curran and van Rees \cite{cvr} that scs(4)=4 and
scs(5)=6, Bate and van Rees \cite{bvr} that scs(6)=9, and Adams
and Khodkar \cite{ak} that scs(7)=12.  Nelder \cite{nelseb}, 
Bate and van Rees \cite{bvr}, and Mahmoodian \cite{mah} independently
conjectured that scs$(n) = \lfloor \frac{n^2}{4} \rfloor$.
Here we show that scs(8)=16.

For a set $S$ in a Latin square $L$ we define sets for each row
$i$, column $j$ and element $k$. Let \ $R_i(S) = \{ k |\ (i, j; k)
\in S \}, \ C_j(S) = \{ k |\ (i, j; k) \in S \}$, \ and \ $E_k(S)
= \{ (i,j) |\ (i, j; k) \in S \}$. So $R_i(S)$ ($C_j(S)$) is the
set of elements which appear in row $i$ (column $j$) of $S$ and
$E_k(S)$ is the set of positions where the element $k$ appears in
$S$.

Now define $R(S) = | \{ i | R_i(S) \neq \emptyset \} |$,
           $C(S) = | \{ j | C_j(S) \neq \emptyset \} |$,
       $E(S) = | \{ k | E_k(S) \neq \emptyset \} |$.

We will write that a trade $T$ has $R(T)$ rows, $C(T)$ columns and
$E(T)$ elements.

Based on results by Cavenagh \cite{njc}, and considering the paper
by Horak, Fleischner, and Aldred \cite{haf}, the author decided to
conjecture the following.

\begin{conjecture}
In a subset $S$ of $L$, a Latin square of order $n \geq 3$, $S
\leq \lfloor \frac{n^2}{4} \rfloor$, there exists a a Latin trade
$T \in L$ such that $T \cap S = \emptyset$ and $R(T) \leq 3, C(T)
\leq 3$, or $E(T) \leq 3$.
\end{conjecture}


\section{Integer programming methods}
A survey of defining sets by Donovan et al \cite{def} includes a
description of an algorithm called ``Algorithm B'' and the
$minimising$ technique reproduced verbatim here.  Smallest
defining sets for designs are a concept analogous to smallest
critical sets for Latin squares, and blocks are here analogous to
entries in Latin squares.

\begin{algorithm}
$\mathbf{(B)}$ Find some trades in $D$ and put these in a list
$\mathcal{T}$.  Now perform steps (1) and (2) repeatedly until $S$
has only one completion.

(1) Form the integer programme corresponding to $\mathcal{T}$,
find an optimal solution, and form the set of blocks $S$
corresponding to this solution.

(2) If $S$ does not complete uniquely, then completions not equal
to $D$ define additional trades in $D$, and these are added to
$\mathcal{T}$.

The running time for step (1) can be improved by $minimising$
$\mathcal{T}$.  That is, if $T_a \in \mathcal{T}$ and $T_b \in
\mathcal{T}$, and if $T_a \subsetneq T_b$, then $T_b$ can be
removed from $\mathcal{T}$.
\end{algorithm}

This algorithm is suitable for finding a smallest defining set for
a design, but it is not suitable for finding the size of the
smallest critical set in Latin squares of order 8, as it requires
too much computation time.  In this paper, we are concerned only
with the $size$ of the smallest critical set.

In the paper by Donovan et al it is stated that ``Algorithm B
yields only a single smallest defining set''.  However, it can be
easily extended with a third step, as follows, to find all
smallest defining sets.

(3) If $S$ has unique completion and contains $s$ blocks, add an
additional constraint $C$ to the integer programme which ensures
that less than $s$ of the blocks in $S$ can occur.  Go back to
step (1).



Similarly to Algorithm B, if $\mathcal{T}$ is the complete set of
Latin trades in a Latin square $L$, then the optimal solution to
the following integer programme is the size of the minimal
critical set in $L$.
\newline

\noindent Minimize: $\sum_{x \in L}{C_{x}}$ \\
Subject to: \\
for each $T \in \mathcal{T}, \sum_{x \in T}{C_{x}} \geq 1$ \\
where $C_{x}$ is 1 if $x \in C$ and 0 otherwise.
\newline

The application of Conjecture 1 means that in order to determine
the value of scs($n$) for small values of $n$, we should first
find representatives of each main class (\cite{klt}, \cite{ak}) of Latin square of order
$n$ and for each square, find the Latin trades on 3 rows, columns
and elements and write a 0-1 integer programme as above.  However,
for our purposes it is not necessary to calculate the exact value
of the smallest critical set for every main class.  In general, it
was found that limiting the size of the trades to a value
determined by trial and error and considering trades only with
less than or equal to 3 rows, columns, or elements was effective.
Thus the following algorithm was used for a Latin square $L$ of
order $n$.

\begin{algorithm}
Find all trades with less than or equal to 3 rows, columns, or
elements.  This is achieved by considering all completions of $L
\setminus S$, where $S$ is a set of 3 rows, columns, or elements
of $L$.  When $S$ ranges over all possible sets of 3 rows,
columns, or elements of $L$, then we can find all such trades by
considering the difference between $L$ and each possible
completion.  See Bean \cite{thesis} for completion methods.
Minimize the list of trades as in Algorithm B.

In the following, the initial value of $z$ and the increment are
chosen to minimize running time.

(1) Create an integer programme $IP$ as above using the trades of
size less than $z$. (Adding an extra constraint to ensure that the
sum of the $C_x's$ does not equal or exceed $\lfloor \frac{n^2}{4}
\rfloor$ speeds up the process.) Solve $IP$.

(2) If $IP$ has a solution with size $< \lfloor \frac{n^2}{4}
\rfloor$ then increase $z$ and repeat step 1, until $z = 3n$.
Otherwise stop.

(3) Similarly to the new step (3) above, check the solution $S$ of
size $s < \lfloor \frac{n^2}{4} \rfloor$ for unique completion. If
it does not have unique completion, remove it from future
consideration by adding a constraint which ensures that less than
$s$ of the cells of $S$ are considered in future iterations.
\end{algorithm}

We used CPLEX \cite{cplex} or BonsaiG \cite{bonsaig} to solve the
integer programme.

scs($n$) may be determined by considering all the main classes of
Latin squares of order $n$, as Donovan et al \cite{dcns} showed.
Kolesova, Lam and Thiel \cite{klt} determined the 283,657 main
classes of Latin squares of order 8.  These, and the main classes
of smaller orders, were obtained from McKay's website \cite{bdm}.

For $n=6$, using Algorithm 2 and limiting the size of the trades
to 10 or less, we can determine that scs(6)=9 in approximately 5
seconds on an Athlon 1200Mhz computer.  Twelve main classes must
be examined, thus each main class takes about half a second.  For
$n=7$, by limiting the size of the trades to 11 or less, we can
determine that scs(7)=12 in 34.5 minutes on the same computer.  We
examine 147 main classes, thus each main class takes about 12
seconds.

As might be expected, there is a strong correlation between the
number of intercalates in a Latin square of order $n$ and the time
taken to prove the non-existence in the square of a critical set
of size less than $\lfloor \frac{n^2}{4} \rfloor$.  For example,
there are two main classes of $7 \times 7$ Latin squares with no
intercalates; in the computation mentioned above, one of these
main classes, $\mathbb{Z}_7$, accounts for 13.8 minutes of the
total 34.5 minute computation. The other main class contains 12
Latin trades of order 6 (the size of the smallest trade in
$\mathbb{Z}_7$ is 9) and accounts for 2.2 minutes. The one main
class with 1 intercalate takes 1.1 minutes, and every other main
class takes 45 seconds or less. In the $6 \times 6$ Latin squares,
the effect is less apparent, but the four main classes with less
than $\lfloor \frac{6^2}{4} \rfloor$ intercalates are the four
Latin squares which take the most time.

For a given order 8 Latin square, it is much easier to test
whether a critical set of size less than 16 exists than it is to
find the size of the smallest critical set in that square.  This
is strictly a constraint satisfaction problem rather than an
integer programming problem, because the minimization step is not
necessary.  Increasing the number of constraints, for this
problem, seems to in general slow solution down.  We began by
limiting the size of the trades to $10$, which was sufficient for
more than 95\% of Latin squares of order 8, before increasing the
limit to 12, 14 and then 24.

For the first 138,800 main classes (arranged in descending order
of the number of intercalates) the computers used were an Athlon
1200Mhz, and computers with 5 Athlon MP 1667Mhz CPUs at IPM,
Tehran.  CPLEX was used. After this a network at the University of
Queensland consisting of 128 Sun Fireblades and 33 dual Pentium
III 800 MHz computers were used with BonsaiG to complete the
calculation, with completion of some gaps by a Pentium IV 2.4 Ghz
at IPM.  The programmes written can be obtained by emailing the
author.  The first part of computation took about three months,
and the second part about one month.

As Conjecture 1 was true for $n=7$ it was conjectured that it
would be true for all $n$.  However, while determining the main
result of this paper, the author discovered the following three
main classes $X$, $Y$, and $W$ (see Table 1).  In $X$ there exist
4 sets $X_1, X_2, X_3,$ and $X_4$, each of size 15, such that for
each Latin trade $T$ in $X_i, 1 \leq i \leq 4$ such that $T \cap
X_i = \emptyset$, $R(T), C(T)$ and $E(T)$ are all greater than 3.
In other words, for each of these sets, every trade in the square
with less than or equal to three rows, columns, or elements
intersects this set.  Similarly, in $Y$ there exist 12 sets of
size 15 $Y_1, \dots, Y_{12}$ with this property, and in $W$ the
set of size 15, $W_1$. These are the only three $8 \times 8$ main
classes with this property, thus step (3) of Algorithm 2 is used
only for these main classes. These sets are presented in Tables 2,
3 and 4.

Thus even if a method could be found to generalize the trades
constructed by Cavenagh to all trades on three rows, columns, or
elements, such a method would not be enough to prove that scs($n$)
$= \lfloor \frac{n^2}{4} \rfloor$.

\section{Properties of the exceptional main classes}
\begin{table}[H]
\begin{center}
\begin{tabular}{ccc}
$
\begin{latinsq}
\hline 1 & 6 & 3 & 4 & 5 & 2 & 8 & 7 \\
\hline 5 & 7 & 1 & 2 & 3 & 8 & 6 & 4 \\
\hline 4 & 5 & 2 & 8 & 6 & 7 & 3 & 1 \\
\hline 3 & 2 & 8 & 7 & 1 & 6 & 4 & 5 \\
\hline 2 & 8 & 4 & 6 & 7 & 1 & 5 & 3 \\
\hline 7 & 4 & 6 & 3 & 8 & 5 & 1 & 2 \\
\hline 8 & 1 & 7 & 5 & 4 & 3 & 2 & 6 \\
\hline 6 & 3 & 5 & 1 & 2 & 4 & 7 & 8 \\
\hline
\end{latinsq}
 $ & $
\begin{latinsq}
\hline 1 & 6 & 3 & 4 & 5 & 2 & 8 & 7 \\
\hline 3 & 5 & 2 & 7 & 4 & 8 & 6 & 1 \\
\hline 4 & 3 & 7 & 8 & 6 & 2 & 1 & 5 \\
\hline 2 & 1 & 8 & 6 & 7 & 3 & 5 & 4 \\
\hline 5 & 7 & 1 & 2 & 8 & 6 & 4 & 3 \\
\hline 7 & 8 & 6 & 1 & 5 & 4 & 3 & 2 \\
\hline 8 & 4 & 5 & 3 & 1 & 7 & 2 & 6 \\
\hline 6 & 2 & 4 & 5 & 3 & 1 & 7 & 8 \\
\hline
\end{latinsq}
$ & $
\begin{latinsq}
\hline 8&6&7&4&1&2&3&5 \\
\hline 7&3&2&8&5&4&1&6 \\
\hline 6&4&5&3&7&1&2&8 \\
\hline 2&5&3&1&8&6&4&7 \\
\hline 5&8&1&2&6&3&7&4 \\
\hline 1&7&4&6&3&8&5&2 \\
\hline 3&2&8&5&4&7&6&1 \\
\hline 4&1&6&7&2&5&8&3 \\
\hline
\end{latinsq}$ \\
\vspace{2mm} \\
X & Y & W \\
\end{tabular}
\caption{The three main classes $X, Y,$ and $W$}
\end{center}
\end{table}

\begin{table}[H]
\begin{center}
\begin{tabular}{cc}
$
\begin{latinsq}
\hline   & 6 & 3 &   &   &   &   &  \\
\hline 5 &   & 1 &   &   &   &   &  \\
\hline   &   &   &   &   & 7 & 3 &  \\
\hline   &   &   & 7 & 1 &   &   &  \\
\hline 2 &   &   & 6 &   &   &   &  \\
\hline   & 4 &   &   &   & 5 &   &  \\
\hline   &   &   &   & 4 &   & 2 &  \\
\hline   &   &   &   &   &   &   & 8\\
\hline
\end{latinsq}
$ & $
\begin{latinsq}
\hline   &   & 3 &   & 5 &   &   &  \\
\hline 5 &   &   &   &   &   & 6 &  \\
\hline 4 &   &   &   &   & 7 &   &  \\
\hline   & 2 &   &   & 1 &   &   &  \\
\hline   &   &   & 6 &   & 1 &   &  \\
\hline   & 4 &   & 3 &   &   &   &  \\
\hline   &   & 7 &   &   &   & 2 &  \\
\hline   &   &   &   &   &   &   & 8\\
\hline
\end{latinsq}
$
\end{tabular}
\\
\vspace{2mm}
\begin{tabular}{cc}
$
\begin{latinsq}
\hline   &   & 3 &   &   & 2 &   &  \\
\hline 5 & 7 &   &   &   &   &   &  \\
\hline   &   &   &   & 6 & 7 &   &  \\
\hline 3 &   &   &   & 1 &   &   &  \\
\hline   &   & 4 & 6 &   &   &   &  \\
\hline   & 4 &   &   &   &   & 1 &  \\
\hline   &   &   & 5 &   &   & 2 &  \\
\hline   &   &   &   &   &   &   & 8\\
\hline
\end{latinsq}
$ & $
\begin{latinsq}
\hline 1 &   &   & 4 &   &   &   &  \\
\hline   &   &   & 2 & 3 &   &   &  \\
\hline   & 5 & 2 &   &   &   &   &  \\
\hline   &   &   &   &   & 6 & 4 &  \\
\hline   &   &   &   & 7 &   & 5 &  \\
\hline 7 &   & 6 &   &   &   &   &  \\
\hline   & 1 &   &   &   & 3 &   &  \\
\hline   &   &   &   &   &   &   & 8\\
\hline
\end{latinsq}
$
\end{tabular}
\caption{$X_1, X_2, X_3$ and $X_4$}
\end{center}
\end{table}

\begin{table}[H]
\begin{center}
\begin{tabular}{ccc}
$
\begin{latinsq}
\hline   &   & 3 & 4 &   &   &   &     \\
\hline   & 5 &   & 7 &   &   &   &     \\
\hline 4 &   &   &   & 6 &   &   &     \\
\hline 2 & 1 &   &   &   &   &   &     \\
\hline   &   & 1 &   &   & 6 &   &     \\
\hline   &   &   &   & 5 &   & 3 &     \\
\hline   &   &   &   &   & 7 & 2 &     \\
\hline   &   &   &   &   &   &   & 8   \\
\hline
\end{latinsq}
$ & $
\begin{latinsq}
\hline   & 6 & 3 &   &   &   &   &     \\
\hline 3 &   &   & 7 &   &   &   &     \\
\hline 4 &   &   &   &   &   & 1 &     \\
\hline   & 1 &   &   & 7 &   &   &     \\
\hline   &   &   & 2 &   & 6 &   &     \\
\hline   &   &   &   & 5 & 4 &   &     \\
\hline   &   & 5 &   &   &   & 2 &     \\
\hline   &   &   &   &   &   &   & 8   \\
\hline
\end{latinsq}
$ & $
\begin{latinsq}
\hline   & 6 &   &   & 2 &   &   &     \\
\hline 3 &   &   &   &   &   & 6 &     \\
\hline   &   & 7 &   &   &   & 1 &     \\
\hline   &   &   &   & 7 & 3 &   &     \\
\hline 5 &   &   & 2 &   &   &   &     \\
\hline   &   &   & 1 &   & 4 &   &     \\
\hline   & 4 & 5 &   &   &   &   &     \\
\hline   &   &   &   &   &   &   & 8   \\
\hline
\end{latinsq}
$
\end{tabular}
\\
\vspace{2mm}
\begin{tabular}{ccc}
$
\begin{latinsq}
\hline 1 &   &   &   & 2 &   &   &     \\
\hline   &   &   &   & 4 &   & 6 &     \\
\hline   &   & 7 &   &   & 2 &   &     \\
\hline   &   &   &   &   & 3 & 5 &     \\
\hline 5 & 7 &   &   &   &   &   &     \\
\hline   &   & 6 & 1 &   &   &   &     \\
\hline   & 4 &   & 3 &   &   &   &     \\
\hline   &   &   &   &   &   &   & 8   \\
\hline
\end{latinsq}
$ & $
\begin{latinsq}
\hline 1 &   &   & 4 &   &   &   &     \\
\hline   & 5 &   &   & 4 &   &   &     \\
\hline   &   &   &   & 6 & 2 &   &     \\
\hline 2 &   &   &   &   &   & 5 &     \\
\hline   & 7 & 1 &   &   &   &   &     \\
\hline   &   & 6 &   &   &   & 3 &     \\
\hline   &   &   & 3 &   & 7 &   &     \\
\hline   &   &   &   &   &   &   & 8   \\
\hline
\end{latinsq}
$ & $
\begin{latinsq}
\hline 1 & 6 &   &   &   &   &   &     \\
\hline 3 &   &   &   & 4 &   &   &     \\
\hline   &   &   &   &   & 2 & 1 &     \\
\hline   &   &   &   & 7 &   & 5 &     \\
\hline   & 7 &   & 2 &   &   &   &     \\
\hline   &   & 6 &   &   & 4 &   &     \\
\hline   &   & 5 & 3 &   &   &   &     \\
\hline   &   &   &   &   &   &   & 8   \\
\hline
\end{latinsq}
$
\end{tabular}
\\
\vspace{2mm}
\begin{tabular}{ccc}
$
\begin{latinsq}
\hline   &   &   & 4 &   & 5 &   &     \\
\hline   & 5 & 2 &   &   &   &   &     \\
\hline   & 3 &   &   & 6 &   &   &     \\
\hline 2 &   &   & 6 &   &   &   &     \\
\hline   &   & 1 &   &   &   & 4 &     \\
\hline 7 &   &   &   &   &   & 3 &     \\
\hline   &   &   &   & 1 & 7 &   &     \\
\hline   &   &   &   &   &   &   & 8   \\
\hline
\end{latinsq}
$ & $
\begin{latinsq}
\hline 1 &   &   &   &   & 5 &   &     \\
\hline   &   & 2 &   & 4 &   &   &     \\
\hline   & 3 &   &   &   & 2 &   &     \\
\hline   &   &   & 6 &   &   & 5 &     \\
\hline   & 7 &   &   &   &   & 4 &     \\
\hline 7 &   & 6 &   &   &   &   &     \\
\hline   &   &   & 3 & 1 &   &   &     \\
\hline   &   &   &   &   &   &   & 8   \\
\hline
\end{latinsq}
$ & $
\begin{latinsq}
\hline   &   & 3 &   & 2 &   &   &     \\
\hline   &   &   & 7 &   &   & 6 &     \\
\hline 4 &   & 7 &   &   &   &   &     \\
\hline   & 1 &   &   &   & 3 &   &     \\
\hline 5 &   &   &   &   & 6 &   &     \\
\hline   &   &   & 1 & 5 &   &   &     \\
\hline   & 4 &   &   &   &   & 2 &     \\
\hline   &   &   &   &   &   &   & 8   \\
\hline
\end{latinsq}
$
\end{tabular}
\\
\vspace{2mm}
\begin{tabular}{ccc}
$
\begin{latinsq}
\hline   &   & 3 &   &   & 5 &   &     \\
\hline   &   & 2 & 7 &   &   &   &     \\
\hline 4 & 3 &   &   &   &   &   &     \\
\hline   & 1 &   & 6 &   &   &   &     \\
\hline   &   &   &   &   & 6 & 4 &     \\
\hline 7 &   &   &   & 5 &   &   &     \\
\hline   &   &   &   & 1 &   & 2 &     \\
\hline   &   &   &   &   &   &   & 8   \\
\hline
\end{latinsq}
$ & $
\begin{latinsq}
\hline   & 6 &   &   &   & 5 &   &     \\
\hline 3 &   & 2 &   &   &   &   &     \\
\hline   & 3 &   &   &   &   & 1 &     \\
\hline   &   &   & 6 & 7 &   &   &     \\
\hline   &   &   & 2 &   &   & 4 &     \\
\hline 7 &   &   &   &   & 4 &   &     \\
\hline   &   & 5 &   & 1 &   &   &     \\
\hline   &   &   &   &   &   &   & 8   \\
\hline
\end{latinsq}
$ & $
\begin{latinsq}
\hline   &   &   & 4 & 2 &   &   &     \\
\hline   & 5 &   &   &   &   & 6 &     \\
\hline   &   & 7 &   & 6 &   &   &     \\
\hline 2 &   &   &   &   & 3 &   &     \\
\hline 5 &   & 1 &   &   &   &   &     \\
\hline   &   &   & 1 &   &   & 3 &     \\
\hline   & 4 &   &   &   & 7 &   &     \\
\hline   &   &   &   &   &   &   & 8   \\
\hline
\end{latinsq}
$
\end{tabular}
\caption{$Y_1,\dots,Y_{12}$}
\end{center}
\end{table}

\begin{table}[H]
\begin{center}
$
\begin{latinsq}
\hline 8 &   & 7 &   &   &   &   &     \\
\hline   & 3 &   &   &   &   & 1 &     \\
\hline 6 &   &   & 3 &   &   &   &     \\
\hline   &   &   &   & 8 &   &   &     \\
\hline   &   &   & 2 &   &   &   & 4   \\
\hline   &   &   &   &   &   & 5 & 2   \\
\hline   &   &   &   & 4 & 7 &   &     \\
\hline   & 1 & 6 &   &   &   &   &     \\
\hline
\end{latinsq}
$
\caption{$W_1$}
\end{center}
\end{table}

We now examine properties of the main class representatives $X$,
$Y$ and $W$.

$X$ and $Y$ are in the class of Latin squares Meynert \cite{am}
refers to as $rcs-symmetric$.  Thus, $(i,j;k) \in L$ implies entry
$(j,k;i) \in L$ and entry $(k,i;j) \in L$.  Although they are not
presented in this way here, they appear in this class.

Each of the sets $X_i$ and $Y_i$ has 11,662,776=$2^3.3^2.161,983$
completions. $W_1$ has 7,075,188 = $2^2.3^5.29.251$ completions.

Each of the sets $X_{i}, Y_{i},$ and $W_{i}$ have $R(j), C(j),$
and $E(j)$ for $j=1,\dots,8$ equal to 1 or 2.  Since the size of
each set is 15, exactly one row and column contains one entry, and
one element occurs once.

$X$ contains 21 intercalates, $Y$ contains 7 intercalates, and $W$
contains 9 intercalates.  Note that $X$ has three times as many
intercalates as $Y$ and $Y$ has three times as many sets as $X$
with the property above.

The sets $Y_{i}$ are symmetric about the main diagonal.

For both $X$ and $Y$, there are seven intercalates intersecting
the entry $(8,8;8)$ and in both cases the sets $X_{i}$ and $Y_{i}$
use in total exactly 43 entries from the squares $X$ and $Y$
respectively.  In $X_{i}$ each entry present occurs once, except
$(8,8;8)$ which occurs 4 times and a transversal (a set in which
every element, row, and column occurs exactly once) in which each
entry (except $(8,8;8)$) occurs 3 times; in $Y_{i}$ each entry
present occurs 4 times, except $(8,8;8)$ which occurs 12 times.

$X$ and $Y$ seem to have a similar spectrum of possible critical
set sizes - in both cases, the spectrum found so far is 21 to 28.
This leads naturally to a question, given that they share so many
other properties.

$\mathbf{Question.}$  Do there exist two main classes of Latin
squares of the same order such that there is a mapping between
critical sets from one main class to the other?

\section{Ideas for future research}
Further research could include looking at the 1707 main classes of
order 9 which contain no intercalates, to check whether any
critical sets of size less than $\lfloor \frac{9^2}{4} \rfloor$
exist there. Alternatively, the computation could be repeated to
verify ``Conjecture 1'' of Bate and van Rees \cite{bvr2} that the
critical set of size $\lfloor \frac{n^2}{4} \rfloor$ exists only
in the main class $\mathbb{Z}_n$, for $n=8$.

\section{Acknowledgements}
The author wishes to thank Auburn University, Discrete and
Statistical Sciences, for their hospitality and the Department of
Mathematics at the University of Queensland for the use of their
computer network.

\end{document}